\documentclass[11pt,reqno]{amsart}
\usepackage{amsmath, amssymb, amsfonts, latexsym}

\newtheorem{prop}{Proposition}
\newtheorem{theo}{Theorem}

\newcommand{\si}{\sigma}
\newcommand{\na}{\nabla}
\newcommand{\om}{\omega}

\newcommand{\La}{\Lambda}

\newcommand{\ka}{K{\"a}hler }

\newtheorem{conjecture*}{Conjecture}
\newtheorem{theorem*}{Theorem}
\newtheorem{question*}{Question}

\newcommand{\leftr}{[\hbox{\hspace{-0.15em}}[}
\newcommand{\rightr}{]\hbox{\hspace{-0.15em}}]}

\newcommand{\bx}{{\bf Q.E.D.}}

\title[Symplectic obstructions to
$\om$-compatible Einstein metrics]{Symplectic obstructions to the
existence of $\om$-compatible Einstein metrics}
\author{TEDI DR\u{A}GHICI}

\address{Department of Mathematics \\ Florida
International University \\ Miami FL 33199 \\ USA}
\email{draghici@fiu.edu}

\begin{document}

\begin{abstract}

It is shown that the existence of an $\om$-compatible Einstein
metric on a compact symplectic manifold $(M,\om)$ imposes certain
restrictions on the symplectic Chern numbers. Examples of
symplectic manifolds which do not satisfy these restrictions are
given. The results offer partial support to a conjecture of
Goldberg.

\vspace{0.1cm} \noindent 2000 {\it Mathematics Subject
Classification}: 53C25
\end{abstract}

\maketitle

\section{Introduction}

This note is motivated by the following still open conjecture of
Goldberg:

\noindent {\bf Conjecture (\cite{Go}):} On a compact symplectic
manifold $(M^{2n}, \om)$ any Einstein $\om$-compatible metric is
K\"ahler Einstein.

A Riemannian metric $g$ is said to be compatible with a
symplectic form $\om$, or shortly,  $\om$-{\it compatible}, if
there exists a $g$-orthogonal almost complex structure $J$ such
that
$$ \om(\cdot, \cdot) = g(J\cdot, \cdot). $$
Such a triple $(g,J,\om)$ is called an {\it almost K\"ahler
structure}.

Given a symplectic form $\om$ on a compact manifold $M^{2n}$, the
space of almost \ka metrics compatible with $\om$ is well known
to be infinite dimensional and contractible. The latter fact
implies that the Chern classes $c_k \in {\bf H}^{2k}(M,
\mathbb{R})$ are independent of the choice of a compatible almost
complex structure. As $\om$ induces a non-trivial cohomology
class $[\om] \in {\bf H}^{2}(M, \mathbb{R})$, we define numerical
symplectic invariants, which we call {\it symplectic Chern
numbers}, by taking cup products of the Chern classes $c_k$ with
appropriate powers of $[\om]$. The symplectic Chern numbers $(c_1
\vee [\om]^{n-1})(M)$ and $(c_1^2 \vee [\om]^{n-2})(M)$ will play
an important role in this note.

It is now well known that \ka metrics exist only very rarely on
compact symplectic manifolds. Indirectly, the Goldberg conjecture
predicts that $\om$-compatible Einstein metrics are even more
rare. Although the conjecture is still wide open, this prediction
can be confirmed in certain cases and our purpose is to bring
further support to its validity.

First, let us mention that for compact 4-manifolds there are
known topological obstructions to the existence of Einstein
metrics. For instance, the Hitchin-Thorpe inequality
$3|\sigma(M)| \leq 2\chi(M)$ should be satisfied, where
$\sigma(M)$, $\chi(M)$ are the signature, respectively, the Euler
number of $M^4$. Important refinements of this inequality were
proved by LeBrun \cite{lebrun1, lebrun3}, using Seiberg-Witten
theory. There are now known many examples of compact symplectic
manifolds which violate the Hitchin-Thorpe inequality or its
refinements and, hence, do not admit {\it any} Einstein metrics
(compatible or not). This provides indirect support to the
4-dimensional Goldberg conjecture. In higher dimensions there are
no known topological obstructions to the existence of Einstein
metrics.

There are results directly supporting the Goldberg conjecture.
Most notably, Sekigawa proved in \cite{sekigawa} that the
conjecture is true provided that the scalar curvature is assumed
to be non-negative. Other positive partial results have been
obtained in dimension 4 under various additional curvature
assumptions \cite{oguro-sekigawa1,oguro-sekigawa,Arm1,AA}.
However these partial results do not provide obstructions to the
existence of Einstein compatible metrics, because of the
Riemannian nature of the additional assumptions imposed.

It was observed in \cite{Dr1}, that Sekigawa's result can be
slightly improved by replacing the assumption $s \geq 0$ with the
weaker condition $(c_1 \vee [\om]^{n-1})(M) \geq 0$. As we need
its proof later on, we incorporate this remark as part of our
main result. Furthermore, in dimension 4, Armstrong proved that
the integrability holds even when one replaces the symplectic
condition $(c_1 \vee [\om])(M) \geq 0$, with, the essentially
topological one, that the manifold admits a metric of everywhere
positive scalar curvature (see \cite{arm}, Corollary 2.3.5).

The main goal of this note is to investigate the case $(c_1 \vee
[\om]^{n-1})(M) < 0$. We prove that the existence of an Einstein
$\om$-compatible metric imposes certain inequalities between the
symplectic Chern numbers $(c_1 \vee [\om]^{n-1})(M)$ and
$(c_1^2\vee [\om]^{n-2})(M)$, which are not satisfied by all
symplectic manifolds. The following theorem summarizes our main
results:
\begin{theo}\label{mainthm} Let $(M^{2n}, \om)$ be a $2n$-dimensional
compact symplectic manifold. Assume that $M$ admits an
$\om$-compatible Einstein metric $g$.
\begin{itemize}
\item{\bf A.} If $(c_1 \vee [\om]^{n-1})(M) \geq 0$, then $g$ is a
K\"ahler-Einstein metric. In particular, $c_1 \in \mathbb{R}_+
[\om]$.
\item{\bf B.} If $(c_1 \vee [\om]^{n-1})(M) < 0$, then the
following inequalities hold:
\begin{equation} \label{ineq1}
(c_1^2 \vee [\omega]^{n-2}(M)) \cdot ([\omega]^{n}(M)) < k_1 ( c_1
\vee [\omega]^{n-1}(M))^2 ,
\end{equation}
where $k_1= 25/9$ if $2n \geq 6$ and $k_1=9/4$ if $2n = 4$;
\begin{equation} \label{ineq2}
(c_1^2 \vee [\omega]^{n-2}(M)) \cdot ([\omega]^{n}(M)) > k_2
 ( c_1 \vee [\omega]^{n-1}(M))^2 ,
\end{equation}
where $k_2= \frac{n-(25/9)}{n-1}$ if $2n \geq 6$ and $k_2=2/3$ if
$2n = 4$.
\end{itemize}
\end{theo}
Part {\bf A} of Theorem 1 leads to first examples of compact
symplectic manifolds of any dimension which do not admit
compatible Einstein metrics. Indeed, any symplectic manifold
$(M,\om)$ with $(c_1 \vee [\om]^{n-1})(M) \geq 0$, but $c_1
\notin \mathbb{R}[\om]$ has this property. Concerning part ${\bf
B}$, the constants $k_1, k_2$ are most likely not optimal. In
fact, I recently learned from Claude LeBrun \cite{lebrun-pr} that
in dimension 4 inequality (\ref{ineq2}) still holds for $k_2 =
3/4$. One would hope the result to be valid with $k_1, k_2$ as
close to 1 as possible. Nevertheless, even with the current
constants, in Section 4 we give examples of symplectic manifolds
which violate (\ref{ineq1}) or (\ref{ineq2}) and thus cannot
admit compatible Einstein metrics.

\section{Preliminaries}

Assume for the beginning that $(M^{2n},g, J, \om)$ is only an
almost Hermitian manifold, i.e. that the fundamental form $\om$ is
not necessarily closed. We shall use the following notations:
$\nabla$ is the Levi-Civita connection, $R$, $Ric$, $s$ are
respectively the curvature tensor, the Ricci tensor and the
scalar curvature of $\nabla$; $\si = \frac{\om^n}{n!}$ is the
volume form and $( \;, \;)$ is the pointwise inner product
induced by the metric $g$ on various bundles of tensors and forms.

The almost complex structure $J$ induces an involution on the
bundle of real 2-forms, by
$$ \Lambda^2 M \ni \xi(\cdot,\cdot) \longrightarrow
\xi(J\cdot,J\cdot) \in \Lambda^2 M .$$ The $\pm 1$-eigenspaces of
this involution, which we denote by $\La^{1,1}_\mathbb{R} M$ and
$\leftr \La^{0,2}M \rightr $, are the bundles of $J$-invariant,
respectively, $J$-anti-invariant 2-forms. The notation is
explained by the correspondence with the usual type decomposition
of complex 2-forms: $J$-invariant 2-forms are nothing but real
forms of complex type (1,1), while $J$-anti-invariant 2-forms are
real parts of complex 2-forms of type (0,2) (equivalently, of
type (2,0)). The fundamental form $\om$ is $J$-invariant and we
denote by $\La^{1,1}_0 M \subset \La^{1,1}_\mathbb{R} M $ the
sub-bundle of {\it primitive} real (1,1)-forms, i.e.
$J$-invariant 2-forms which are point-wise orthogonal to $\om$.
Thus we have
\begin{equation}\label{Lambda2r}
\La^2M = \La^{1,1}_\mathbb{R} M \oplus \leftr \La^{0,2}M \rightr =
(\mathbb{R} \om \oplus \La^{1,1}_0 M) \oplus \leftr \La^{0,2}M
\rightr ,
\end{equation}
and the components of a section $\xi \in \La^2 M$ with respect to
this decomposition are
$$ \xi = \xi' + \xi'' = \frac{1}{n} (\xi, \om) \om + \xi_0' + \xi'' .$$
Here and throughout the paper we use the superscripts $'$ and
$''$ to denote respectively the $J$-invariant and
$J$-anti-invariant components and the subscript $0$ for the
primitive part.

\vspace{0.2cm}

\noindent For any $\xi \in \La^2 M$, easy computations imply:
\begin{equation} \label{wedge1}
\xi \wedge \om^{n-1} = \frac{1}{n} (\xi, \om) \; \om^{n} = (n-1)!
\; (\xi, \om)\; \si ;
\end{equation}
\begin{equation} \label{wedge2}
\xi \wedge \xi \wedge \om^{n-2} = (n-2)! \; \Big[ \frac{n-1}{n}
(\xi, \om)^2 - |\xi_{0}'|^2 + |\xi''|^2 \Big] \si =
\end{equation}
$$ = (n-2)! \; \big[ (\xi, \om)^2 - |\xi'|^2 + |\xi''|^2 \big] \si .$$


From now on we assume that $(g, J, \omega)$ is an almost \ka
structure, i.e that $\om$ is closed. It is well known that for an
almost \ka structure, $\nabla \om$ is identified with the
Nijenhuis tensor $N$ of $J$ by (cf. e.g. \cite{KN}):
\begin{equation} \label{naom-N}
(\na_X\om)(\cdot,\cdot) = \frac{1}{2} ( JX, N(\cdot,\cdot) ).
\end{equation}
Since $N(J\cdot, \cdot) = N(\cdot, J\cdot)= - JN(\cdot, \cdot)$,
the identification (\ref{naom-N}) implies that for any tangent
vectors $X,Y,Z$
\begin{eqnarray} \label{naom-anti}
(\na_X \om)(JY, JZ) &=& - (\na_X \om)(Y,Z) ;\\
\label{qk} (\na_{JX} \om)(JY,Z) &=& - (\na_X \om)(Y,Z) .
\end{eqnarray}
Relation (\ref{qk}) is sometimes called the {\it quasi-\ka
condition}. The trace in $X,Y$ of (\ref{qk}) leads to the (again
well known) fact that $\om$ is also co-closed and hence harmonic
with respect to $g$.

The standard Weitzenb{\"o}ck formula for 2-forms
\begin{equation*} \label{wtz2}
\Delta \xi -\na^* \na \xi = [{\rm Ric}(\xi \cdot,\cdot) - {\rm
Ric}(\cdot, \xi\cdot)] - 2 {R}(\xi),
\end{equation*}
specialized to $\xi = \om$, gives
\begin{equation} \label{wtz2Om}
\frac{1}{2} \na^* \na \om = {R} (\om) - \frac{1}{2}[{\rm
Ric}(J\cdot, \cdot) - {\rm Ric}(\cdot, J\cdot )] = \rho_* - \rho
\; .
\end{equation}
Formula (\ref{wtz2Om}) is a measure of the difference of two
types of Ricci forms. For an arbitrary almost \ka structure the
Ricci tensor is in general not $J$-invariant, but taking its
$J$-invariant part ${\rm Ric}'$, we can define the {\it Ricci
form}, $\rho(\cdot, \cdot) = {\rm Ric}'(J\cdot, \cdot)$. The
2-form defined by $\rho_* = R(\om)$ is called the {\it $*$-Ricci
form}; this is in general not $J$-invariant. In fact, it follows
from (\ref{wtz2Om}) that $\rho_*'' = \frac{1}{2} (\nabla^* \nabla
\om)''$. As for the $J$-invariant part of (\ref{wtz2Om}), taking
the covariant derivative $\na_W$ of the relation
(\ref{naom-anti}) and then taking the trace in $W, X$, we obtain
$ (\nabla^* \nabla \om)' = \psi$, where $\psi$ is the
semi-positive 2-form given by
$$ \psi(X,Y) = \sum_{i=1}^{2n} ( (\na_{e_i} J)JX, (\na_{e_i} J)Y )
.$$ Here and throughout $\{e_i\}_{i=1, 2n}$ denotes an
orthonormal basis with respect to $g$. A $J$-invariant 2-form $\xi
\in \Lambda^{1,1}_{\mathbb{R}} M$ is called {\it semi-positive}
if $\xi(X,JX) \geq 0, \; \; \forall X \; \in \; TM$.

The inner product with $\om$ of the relation (\ref{wtz2Om})
yields the difference of the two types of scalar curvatures:
\begin{equation} \label{s*s}
s^* - s = |\na \om|^2 = \frac{1}{2}|\na J|^2 \; ,
\end{equation}
where $s^*=2(R(\om), \om)$, is the so called $*$-{\it scalar
curvature}.

\vspace{0.2cm}

Unlike the \ka case, the Levi-Civita connection cannot be used
directly to provide representatives for the Chern classes $c_k$.
Instead, one uses the so called {\it Hermitian} or {\it first
canonical} connection (see e.g. \cite{gauduchon1}), defined by :
$${\widetilde{\na}}_X Y = \na_X Y - \frac{1}{2}J(\na_X J)(Y) . $$
If $\widetilde{R}$ denotes the curvature tensor of ${\widetilde
{\na}}$, then $$ \widetilde{\rho}(X,Y) = \frac{1}{2}
\sum_{i=1}^{2n}(\widetilde{R}_{X,Y} e_i, Je_i)$$ is a closed
2-form which is a deRham representative of $2\pi c_1$ in
$H^2(M,{\mathbb R})$. One easily finds the explicit relation
between the curvature tensors $\widetilde{R}$ and $R$, of
${\widetilde {\na}}$ and $\na$. We will only need the
relationship of the Ricci forms:
\begin{equation} \label{tilde}
\widetilde{\rho} = \rho^{*} - \frac{1}{2} \phi ,
\end{equation}
where $\phi$ is the $J$-invariant, semi-positive 2-form given by
$\phi(X,Y) = (\na_{JX} \om, \na_Y \om)$.

\vspace{0.2cm}

 Hence, by (\ref{wedge1}), (\ref{wedge2}),
(\ref{s*s}) and (\ref{tilde}) we have
\begin{equation}  \label{c1om}
\frac{4\pi}{(n-1)!} (c_1 \vee [\omega]^{n-1})(M) = \int_M
\frac{1}{2}(s^* + s) \; \si=\int_M (s + \frac{1}{2} |\na \om|^2)
\; \si,
\end{equation}
\begin{equation}  \label{c1^2om}
\frac{4\pi^2}{(n-2)!} (c_1^2 \vee [\omega]^{n-2})(M) = \int_M
\Big[ \frac{(s^* + s)^2}{16} - |\rho_*' - \frac{1}{2} \phi|^2 +
|\rho_*''|^2 \Big] \; \si .
\end{equation}
The formula (\ref{c1om}) is due to Blair \cite{blair}, who first
noted that the integral $\int_M (s^* + s) \; \si$ is a symplectic
invariant. We let the reader observe that formulas (\ref{c1om})
and (\ref{c1^2om}) reduce to the well known ones in the K\"ahler
case.

\vspace{0.2cm}

\noindent We close this section with the following classical
result of Apte about the Chern numbers $(c_1^2 \vee
[\om^{n-2}])(M)$ and $(c_1 \vee [\om^{n-1}])(M)$ in the K\"ahler
case:
\begin{prop} \cite{apte} \label{prop1} Let $M^{2n}$ be a compact
manifold and let $\om$ be a symplectic form on $M$ which admits a
compatible \ka metric. Then
\begin{equation} \label{kaineq}
(c_1^2 \vee [\om]^{n-2})(M) \cdot ([\om]^n)(M) \leq ((c_1 \vee
[\om]^{n-1})(M))^2 ,
\end{equation}
with equality iff $c_1 \in \mathbb{R} [\om]$.
\end{prop}

To  sketch a proof (slightly different than the original one of
\cite{apte}, see also \cite{Pe}), note that in the \ka case the
decomposition (\ref{Lambda2r}) descends to cohomology. In view of
(\ref{wedge2}), the bilinear form $b(c,d) = (c \vee d \vee
[\om]^{n-2})(M)$ has Lorenz signature $(+,-,...,-)$ when
restricted to ${\bf H}_{\mathbb{R}}^{1,1} \times {\bf
H}_{\mathbb{R}}^{1,1}$, where ${\bf H}_{\mathbb{R}}^{1,1}$
denotes the subset of ${\bf H}^2(M, \mathbb{R})$ consisting of
cohomology classes represented by real harmonic 2-forms of type
(1,1). This fact is part of the so-called Hodge-Riemann bilinear
relations (see \cite{gr-ha}, p. 123). For any $c \in {\bf
H}_{\mathbb{R}}^{1,1}$, we then have the following ``opposite''
Cauchy-Schwarz inequality:
$$ b(c, c) \cdot b([\om] , [\om]) \leq (b(c, [\om]))^2 , $$
with equality iff $c \in \mathbb{R} [\om]$. It is well known that
for a \ka manifold the first Chern class $c_1$ belongs to ${\bf
H}_{\mathbb{R}}^{1,1}$.

The proposition is no longer true in the non-\ka case. One can
find examples of symplectic forms which do not satisfy the
conclusion of the Proposition 1, and hence do not admit
compatible \ka metrics (see \cite{Dr2} and Proposition \ref{exp1}
below).

\section{Proof of Theorem 1}

We start by recalling the remarkable integral formula of Sekigawa,
which is valid on an arbitrary compact almost \ka manifold. The
original proof of this formula \cite{sekigawa} is based on
Chern-Weil theory. An alternative approach, based on Weitzenb\"ock
formulae was described in \cite{ADM}.
\begin{prop}\label{o1}\cite{sekigawa} For
any compact almost-\ka manifold $(M^{2n}, g, J, \om)$, the
following integral formula holds:
\begin{equation} \label{seki1} 0 = \int_M \Big[ \frac{1}{2}|{\rm
Ric}''|^2 - |\rho_*''|^2 - 2|W''|^2 + ({\rho}, \phi - \psi) -
\frac{1}{4}|\psi|^2 - \frac{1}{4}|\phi|^2  \Big] \si \; .
\end{equation}
\end{prop}
The notations are those from section 2; we should add that $W''$
is a certain component of the Weyl part of the curvature (for more
details see \cite{ADM}). For our purposes here, all that matters
is that $|W''|^2$ is a non-negative quantity.

According to (\ref{Lambda2r}),
$$ \rho = \frac{s}{2n} \om + \rho_0 , \; \; \phi =
\frac{|\na \om|^2}{2n} \om + \phi_0 , \; \; \psi = \frac{|\na
\om|^2}{n} \om + \psi_0 ,$$ hence (\ref{seki1}) becomes
\begin{eqnarray} \label{seki1'} 0 &=& \int_M \Big[ \frac{1}{2}|{\rm
Ric}''|^2 - |\rho_*''|^2  - 2|W''|^2 + ({\rho}_0, \phi_0 - \psi_0)
\\ \nonumber & &  - \frac{s}{4n} |\na \om|^2 - \frac{5}{16n} |\na \om|^4
- \frac{1}{4}|\psi_0|^2 - \frac{1}{4}|\phi_0|^2  \Big] \si \; .
\end{eqnarray}
In the Einstein case this implies
$$\int_M ( - s |\na \om|^2) \; \si \geq  \frac{5}{4} \int_M |\na \om|^4  \;
\si $$ and concludes the proof of Sekigawa's theorem
 that compact almost \ka Einstein manifolds with
$s\geq 0$ are necessarily \ka Einstein \cite{sekigawa}. Making no
assumption on the sign of the (constant) scalar curvature and
using Schwarz inequality, one obtains
$$ - s \; {\rm vol}(M) \int_M |\na \om|^2 \; \si
\geq \frac{5}{4} \Big(\int_M |\na \om|^2 \; \si \Big)^2 .$$
Assuming now that the manifold is {\it not} K\"ahler, this leads
to
$$- s \; {\rm vol}(M) \geq \frac{5}{4} \int_M |\na \om|^2 \;
\si ,$$ and, further, using Blair's formula (\ref{c1om}), to
\begin{equation} \label{s-estimate}
(c_1 \vee [\om]^{n-1})(M) > \frac{(n-1)!}{4\pi} \; s \; {\rm
vol}(M) \geq \frac{5}{3} (c_1 \vee [\om]^{n-1})(M).
\end{equation}
In particular, $(c_1 \vee [\om]^{n-1})(M) < 0$, hence part {\bf A}
of Theorem 1 follows by contra-position.

The constant $5/3$ in (\ref{s-estimate}) can be lowered in the
4-dimensional case. In this dimension, the bundle of 2-forms also
decomposes $\La^2 M = \La^+M \oplus \La^-M$, into the sub-bundles
of self-dual and anti-self-dual 2-forms. This is related to the
type decomposition (\ref{Lambda2r}) by
$$\La^+M = \mathbb{R} \om
\oplus \leftr \La^{0,2}M \rightr, \; \; \La^-M = \La^{1,1}_0 M. $$
One then immediately concludes that $(\na^* \na \om )'$ must be a
multiple of $\om$. Also, using (\ref{qk}) and the fact that the
sub-bundle $\leftr \La^{0,2}M \rightr$ has dimension 2, it
follows that the symmetric 2-tensor $(\na_{\cdot} \om,
\na_{\cdot} \om)$ has a double eigenvalue $0$ and a double
eigenvalue $\frac{|\na \om|^2}{2}$. Hence, in dimension 4 we have
\begin{equation} \label{psi0phi0}
\psi_0 = 0, \; \; |\phi_0|^2 = \frac{1}{8} |\na \om|^4.
\end{equation}
Using these in (\ref{seki1'}) and following the path described
above, we obtain that a 4-dimensional Einstein strictly almost \ka
manifold satisfies
\begin{equation} \label{s-estim-4dim}
(c_1 \vee [\om])(M) > \frac{1}{4\pi} \; s \; {\rm vol}(M) \geq
\frac{3}{2} (c_1 \vee [\om])(M).
\end{equation}

\vspace{0.2cm}

Part {\bf B} of Theorem 1 is a consequence of the following
proposition, which may be of interest in its own.
\begin{prop} Let $(M^{2n}, g, J, \om)$ be a compact almost-\ka manifold.
Then the following lower estimates of the $L^2$-norm of the Ricci
tensor hold:
\begin{eqnarray}  \label{L2Ricci1}
\int_M |{\rm Ric}|^2 \; \si &\geq& \frac{8\pi^2}{(n-1)!}\Big(
\frac{n(c_1 \vee [\omega]^{n-1}(M))^2}{[\om]^n(M)} - (n-1)(c_1^2
\vee [\omega]^{n-2})(M) \Big) \\
\label{L2Ricci2} \int_M |{\rm Ric}|^2 \; \si &\geq&
\frac{8\pi^2}{(n-1)!} \; (c_1^2 \vee [\omega]^{n-2})(M).
\end{eqnarray}
Equality holds in (\ref{L2Ricci1}) if and only if $(g,J,\om)$ is
\ka with constant scalar curvature and equality holds in
(\ref{L2Ricci2}) if and only if $(g,J,\om)$ is K\"ahler Einstein.
\end{prop}
\noindent {\it Proof:} Note first that using (\ref{wtz2Om}), we
have
\begin{eqnarray} \nonumber
|\rho_*^{'} - \frac{1}{2} \phi|^2 &=& |\rho + \frac{1}{2}(\psi -
\phi)|^2 = \\ \nonumber &=& |\rho|^2 - <\rho, \phi> + <\rho,\psi>
+ \frac{1}{4}|\psi - \phi|^2 .
\end{eqnarray}
With this, Sekigawa's formula (\ref{seki1}) can also be written as
\begin{equation} \label{seki2} 0 = \int_M \Big[ \frac{1}{2}|{\rm
Ric}|^2  - 2|W''|^2 - \frac{1}{2}<\psi, \phi> - |\rho_*''|^2  -
|\rho_*' - \frac{1}{2} \phi|^2 \Big] \si \; .
\end{equation}
Using (\ref{seki2}) to successively substitute terms in
(\ref{c1^2om}), we get the following alternative expressions for
the symplectic Chern number $(c_1^2 \vee [\om]^{n-2})(M)$:
\begin{eqnarray} \label{c1^2om-1}
\frac{4\pi^2}{(n-2)!} (c_1^2 \vee [\omega]^{n-2})(M) &=& \int_M
\Big[ \frac{(s^* + s)^2}{16} + 2 |\rho_*''|^2 - \frac{1}{2}|{\rm
Ric}|^2  \\ \nonumber & &+ 2|W''|^2 + \frac{1}{2}<\psi, \phi>
\Big] \; \si ; \\ \label{c1^2om-2} \frac{4\pi^2}{(n-2)!} (c_1^2
\vee [\omega]^{n-2})(M) &=& \int_M \Big[ \Big(1-\frac{2}{n}\Big)
\frac{(s^* + s)^2}{16}  +
\frac{1}{2}|{\rm Ric}|^2  \\
\nonumber & & - 2|(\rho_* ' - \frac{1}{2} \phi)_0|^2 - 2|W''|^2 -
\frac{1}{2}<\psi, \phi> \Big] \; \si .
\end{eqnarray}
Since both $\phi$ and $\psi$ are semi-positive 2-forms, relations
(\ref{c1^2om-1}) and (\ref{c1^2om-2}) imply immediately the
following inequalities:
\begin{eqnarray}  \label{ineq1c1^2}
\frac{4\pi^2}{(n-2)!} (c_1^2 \vee [\omega]^{n-2})(M) &\geq& \int_M
\Big[ \frac{(s^* + s)^2}{16} - \frac{1}{2} |{\rm Ric}|^2 \Big] \;
\si \\ \label{ineq2c1^2} \frac{4\pi^2}{(n-2)!} (c_1^2 \vee
[\omega]^{n-2})(M) &\leq& \int_M \Big[ \Big(1-
\frac{2}{n}\Big)\frac{(s^* + s)^2}{16} + \frac{1}{2} |{\rm
Ric}|^2  \\ \nonumber & & - 2 |(\rho_*' - \frac{1}{2} \phi)_0|^2
\Big] \; \si.
\end{eqnarray}
The estimate (\ref{L2Ricci1}) follows from (\ref{ineq1c1^2}),
using (\ref{c1om}) and Schwarz inequality. The estimate
(\ref{L2Ricci2}) follows from  (\ref{ineq2c1^2}) $-
(1-2/n)$(\ref{ineq1c1^2}).

For the equality statement, note first that equality holds in
(\ref{ineq1c1^2}) or (\ref{ineq2c1^2}) if and only if the
structure is K\"ahler. Indeed, assuming equality in either case,
we must have $<\phi, \psi> = 0$. Since both $\phi$ and $\psi$ are
semi-positive, it follows that for any $X \in TM$, $\phi(X,JX) =
0$ or $\psi(X,JX) = 0$. But $\phi(X,JX) = 0$ implies, by the
definition of $\phi$, that $\na_X \om = 0$. The condition
$\psi(X,JX) = 0$ leads to $(\na_Y \om)(X,Z) = - (\na_Y \om)(Z,X)
= 0$, for any $Y,Z \in TM$. But, since $\om$ is closed, this also
leads to $\na_X \om = 0$.

Now further note that for (\ref{L2Ricci1}) we also used Schwarz
inequality, hence in the equality case we must have $s= const$,
while for (\ref{L2Ricci2}) we neglected the last term of
(\ref{ineq2c1^2}), which in the equality case implies $Ric_0 = 0$.
$\bx$

\vspace{0.2cm}

{\it Remark:} Note that the right hand-side of (\ref{L2Ricci1})
is greater or smaller than the right hand-side of
(\ref{L2Ricci2}) depending on whether the inequality
(\ref{kaineq}) holds or not. In the almost \ka case either
situation is possible as it will become clear in Section 4 (see
also \cite{Dr2}).

\vspace{0.2cm}

{\it Proof of Theorem 1, {\bf B.}} In case of dimension $2n \geq
6$, both inequalities (\ref{ineq1}) and (\ref{ineq2}) are now
immediate. Indeed, assuming that $(g,J,\om)$ is a non-K\"ahler,
Einstein, almost \ka structure, by Sekigawa's theorem and Theorem
1, part {\bf A}, both $s$ and $(c_1 \vee [\omega]^{n-1})(M)$ are
negative numbers. The second part of (\ref{s-estimate}) squared
implies then
$$ \frac{((n-1)!)^2}{16\pi^2} \; s^2 \; ({\rm
vol}(M))^2 \leq \frac{25}{9} ((c_1 \vee [\om]^{n-1})(M))^2.$$ Now
combine this inequality with the (strict) inequalities
(\ref{L2Ricci1}) and (\ref{L2Ricci2}) written in the Einstein
case. Inequalities (\ref{ineq1}) and (\ref{ineq2}) follow, with
$k_1=25/9$ and $k_2 =\frac{n - (25/9)}{n-1}$ as stated.

With the same arguments as above, the better constant $k_1= 9/4$
for inequality (\ref{ineq1}) in the 4-dimensional case follows
from (\ref{L2Ricci1}) combined with (\ref{s-estim-4dim}).

To obtain the constant $k_2 = 2/3$ for inequality (\ref{ineq2}) in
dimension 4, slightly more effort is required. With the Einstein
assumption and taking into account (\ref{psi0phi0}), relation
(\ref{c1^2om-1}) becomes
\begin{equation*}
4\pi^2 c_1^2(M) = \int_M \Big[ \frac{(s^* + s)^2}{16}  -
\frac{s^2}{8} + \frac{|\na \om|^4}{8} + 2|\rho_*''|^2 + 2|W''|^2
\Big] \; \si .
\end{equation*}
Using (\ref{s*s}), the above can be written as
\begin{equation*}
4\pi^2 c_1^2(M) = \int_M \Big[ \frac{1}{48}(s^* + s)^2 +
\frac{1}{48} (2s^* - s)^2  + 2|\rho_*''|^2 + 2|W''|^2 \Big] \;
\si ,
\end{equation*}
hence
$$ 4\pi^2 c_1^2(M) > \frac{1}{48} \int_M  (s^* + s)^2  \;
\si . $$ The inequality is strict because if not the structure
would be K\"ahler (see for e.g. \cite{Arm1}), and we assumed
otherwise. Further, using Schwarz inequality and (\ref{c1om}), we
get
$$ (c_1^2(M)) \cdot ([\omega]^{2}(M)) >
\frac{2}{3} ( (c_1 \vee [\omega])(M))^2 , $$ which is the
inequality claimed. $\bx$

\section{Examples}

We already remarked in the introduction that part {\bf A} of
Theorem 1 provides first examples of symplectic manifolds which
do not admit compatible Einstein metrics. We now give such
examples with $(c_1 \vee [\om]^{n-1})(M) < 0$. The first source is
the following proposition, which is essentially inspired from
\cite{Dr2}, but complements the results there.
\begin{prop} \label{exp1} Let $(M^{2n}, J)$ be a compact complex manifold
and assume that $\omega$ is a K\"ahler form  and $\beta$ is a
holomorphic (2,0) form on $(M^{2n},J)$. Then for any $t \in
\mathbb{R}$, the form $ \omega_t = \omega + t Re(\beta)$ is a
symplectic form on $M^{2n}$. Furthermore, if we assume that $c_1
= - [\omega]$, then the following hold:

(i) If $n=2m$ and $\beta^m $ is not identically $0$, then for
$|t|$ large enough, $(M^{4m}, \omega_t)$ does not satisfy
inequalities (\ref{kaineq}) and (\ref{ineq1}), hence it does not
admit compatible K\"ahler metrics, nor compatible Einstein
metrics.

(ii) If $n=2m+1$ and $\beta^m $ is not identically $0$, then for
$|t|$ large enough, $(M^{4m+2}, \omega_t)$ does not satisfy
inequality (\ref{ineq2}), hence it does not admit compatible
Einstein metrics.

(iii) If $n=2m$ or $n=2m+1$ and the highest non-zero power of
$\beta$ is $k < m$, with $(25/9)(n-2k) < n$, then, for $|t|$ large
enough, $(M^{2n}, \omega_t)$ does not satisfy inequality
(\ref{ineq2}), hence it does not admit compatible Einstein
metrics.
\end{prop}

\noindent {\it Proof:} It is well known that on a K\"ahler
manifold any holomorphic form is closed. Thus, $\omega_t$ is
closed for any $t$. To check the non-degeneracy, observe that the
only non-vanishing terms from the binomial expansion of
$\omega_t^n$ are those of the form $\omega^{n-2l}\wedge \beta^l
\wedge \overline{\beta}^l$. But for any form $\alpha$ of type
$(2l,0)$, we have the (pointwise) Hodge-Riemann bilinear relation
(see \cite{gr-ha}, p. 123 and p. 110)
\begin{equation} \label{holoposi}
\omega^{n-2l}\wedge \alpha \wedge \overline{\alpha} = (n-2l)! \,
|\alpha|^2 \frac{\omega^n}{n!},
\end{equation}
where the norm is the one induced by the K\"ahler metric
corresponding to $(J,\omega)$. Thus $\omega_t$ is a symplectic
form for any $t$.

Assuming now that $c_1 = - [\omega]$, the statements from (i),
(ii) and (iii) follow by computing
$$ L= \lim_{t \rightarrow \pm \infty} \frac{ (c_1^2 \vee [\omega_t]^{n-2}(M))
\cdot ([\omega_t]^n(M)) }{ ((c_1 \vee [\omega_t]^{n-1})(M))^2} ,$$
in each case. This is easily accomplished identifying the top
powers of $t$ in the following binomial expansions
$$ \omega_t^n  = \sum_{l=0}^{[\frac{n}{2}]} C_{n}^{2l} C_{2l}^{l} (t/2)^{2l}
\omega^{n-2l} \wedge \beta^l \wedge \overline{\beta}^l ;$$
$$ \omega \wedge \omega_t^{n-1} =
\sum_{l=0}^{[\frac{n-1}{2}]} C_{n-1}^{2l} C_{2l}^{l} (t/2)^{2l}
\omega^{n-2l} \wedge \beta^l \wedge \overline{\beta}^l ;$$
$$ \omega^2 \wedge \omega_t^{n-2} =
\sum_{l=0}^{[\frac{n-2}{2}]} C_{n-2}^{2l} C_{2l}^{l} (t/2)^{2l}
\omega^{n-2l} \wedge \beta^l \wedge \overline{\beta}^l .$$ It
follows that $L= +\infty$ in case (i), $L=0$ in case (ii) and $L=
\frac{n(n-2k-1)}{(n-1)(n-2k)}$ in case (iii). Now the statements
are clear, noting in case (iii) that
$\frac{n(n-2k-1)}{(n-1)(n-2k)} < \frac{n-(25/9)}{n-1}
\Leftrightarrow (25/9)(n-2k) < n$. $\bx$

\vspace{0.2cm}

{\it Remarks:} (a) Certainly the condition $c_1 = -[\omega]$
cannot be replaced by $c_1 = [\omega]$, as in that case there are
no non-trivial holomorphic forms by Kodaira's vanishing theorem.
One would like to understand better the condition that $\beta^m$
is not identically $0$, for $n=2m$, or $n=2m+1$. This is
trivially satisfied if $n=2$, or $3$, by any non-trivial
holomorphic $(2,0)$ form. Further, the condition is stable under
products: if $\beta_1$ has this property on $M_1$ and $\beta_2$
on $M_2$, then so does $\beta_1 + \beta_2$ on $M_1 \times M_2$.
However, for product manifolds (or holomorphic fiber bundles)
case (iii) occurs when $\beta$ is a holomorphic $(2,0)$ form
coming from one of the factors (or from the base).

(b) With the notations from the above proposition, we showed in
\cite{Dr2} that if $(M^4, J, \omega)$ is a compact K\"ahler
surface with $c_1 = -[\omega]$, then for {\it all} values of $t
\neq 0$, the symplectic forms $\omega_t$ violate inequality
(\ref{kaineq}), hence they do not admit compatible K\"ahler
metrics. The same was shown to be true in all higher dimensions
for {\it small} non-zero values of $t$. Now we obtain the same
conclusion when $|t|$ is sufficiently large and $n$ is even. It
is perhaps tempting to conjecture that in any dimension and for
any holomorphic (2,0) form $\beta$, the symplectic 2-form
$\omega_t$ does not admit compatible K\"ahler metrics, for any
$t\neq 0$.

\vspace{0.2cm}

The next source of examples is the following proposition,
suggested to me by Claude LeBrun.

\begin{prop} \label{exp2}
Let $(M^{2n_1}_1, \eta)$,  $(M^{2n_2}_2, \mu)$ be symplectic
manifolds such that
\newline $c_1(M_1) = - [\eta]$, $c_1(M_2) = -
[\mu]$. On $M^{2n} = M^{2n_1}_1 \times M^{2n_2}_2$ ($n=n_1+n_2$),
consider the symplectic forms $\omega_t = \eta + t \mu$, for $t
> 0$. Then the manifold $(M^{2n}, \omega_t)$ does not satisfy
inequality (\ref{ineq2}) and, hence, does not admit compatible
Einstein metrics, in any of the following cases:

(i) if $2n_1= 4$,  and $t$ is sufficiently large;

(ii) if $2n_2 = 4$, and $t$ is sufficiently small;

(iii) if $2n_1 \geq 6, 2n_2 \geq 6$, $(25/9) n_1 < n$, and $t$ is
sufficiently large;

(iv) if $2n_1 \geq 6, 2n_2 \geq 6$, $(25/9) n_2 < n$, and $t$ is
sufficiently small.

\end{prop}

\noindent {\it Proof:} First note that cases (ii) and (iv) can be
obtained from (i), respectively (iii), by substituting $t$ with
$1/t$. For (i) and (iii), we compute as in the previous
proposition
$$ L= \lim_{t \rightarrow \infty} \frac{ (c_1^2 \vee [\omega_t]^{n-2}(M))
\cdot ([\omega_t]^n(M)) }{ ((c_1 \vee [\omega_t]^{n-1})(M))^2} .$$
Note that $c_1(M) = - ([\eta]+ [\mu])$. We easily obtain
$$ \omega_t^n = C_n^{n_1} t^{n_2} \; \eta^{n_1} \wedge \mu^{n_2},$$
$$ (\eta+\mu)^2 \wedge \omega_t^{n-2} =
(C_{n-2}^{n_1-2} t^{n_2} + 2C_{n-2}^{n_1-1} t^{n_2-1} +
C_{n-2}^{n_1} t^{n_2-2}) \; \eta^{n_1} \wedge \mu^{n_2} ,$$
$$ (\eta+\mu) \wedge \omega_t^{n-1} =
(C_{n-1}^{n_1-1} t^{n_2} + C_{n-1}^{n_1} t^{n_2-1} ) \; \eta^{n_1}
\wedge \mu^{n_2} ,$$ with the convention that a binomial
coefficient $C_a^b$ is $0$ if $a\leq 0$, or $b < 0$, or $a< b$.
It follows that $L=0$ in case (i) and $L =
\frac{n(n_1-1)}{n_1(n-1)}$ in case (iii) and the statements are
now clear. $\bx$

\vspace{0.2cm} \noindent {\bf Acknowledgment:} The author is
grateful to Claude LeBrun for suggesting Proposition \ref{exp2}
and for helpful discussions and to Vestislav Apostolov for useful
comments on an earlier version of the paper.

\end{document}